\documentclass{siamart251216}
\usepackage{graphicx} % Required for inserting images
\usepackage{subcaption}
\usepackage{srcltx,epsfig}
\usepackage{amsmath, amssymb}
\usepackage{color}
\usepackage{float}
\usepackage{xcolor}
\usepackage{multirow}
\usepackage{mathtools}
\usepackage{verbatim}
\usepackage{pgfplots}
\pgfplotsset{compat=1.18} 
\usepackage{mathrsfs}
\usepackage{array}
\usepackage{booktabs} 
\usepackage{placeins}

\newtheorem{remark}{Remark}[section]

\usepackage{bm}
\hypersetup{
    colorlinks=true,
    linkcolor=blue,
    anchorcolor=blue,
    citecolor=red
}
\usepackage[T1]{fontenc}
\usepackage{tikz}
\usetikzlibrary{shapes.geometric, arrows}
\usetikzlibrary{matrix}
\usepackage{algorithm}
\usepackage{algpseudocode}
\usepackage{makecell}
\usepackage{siunitx}
\usepackage{chngcntr}
\counterwithin{figure}{section}
\usepackage{indentfirst}
\numberwithin{equation}{section} 

\newcommand{\bmx}{\bm{x}}
\newcommand{\bmX}{\bm{X}}
\newcommand{\bms}{\bm{\sigma}}
\newcommand{\bmb}{\bm{b}}
\newcommand{\bmW}{\bm{W}}
\newcommand{\bmD}{\bm{D}}

\graphicspath{{images/}}

\title{Weak adversarial tensor neural network for high-dimensional Fokker-Planck equations}
\author{
   Sihong Shao\thanks{CAPT, LMAM, and School of Mathematical Sciences, Peking University, Beijing, 100871, China (sihong@math.pku.edu.cn).} \and 
    Zhongwei Xu\thanks{School of Mathematical Sciences, Peking University, Beijing, 100871, China (xuzw@pku.edu.cn).}
    } 
\headers{Weak adversarial tensor neural network}
        {Sihong Shao and Zhongwei Xu}
\begin{document}
\date{}
\maketitle
%\tableofcontents
%\newpage

\begin{abstract}
Solving high-dimensional stationary Fokker–Planck (SFP) equations remains challenging. Recently, deep learning methods have shown great potential for solving high-dimensional SFP equations due to their universal approximation and mesh-free training. However, the high-dimensional integrals involved in loss computation and normalization can be expensive to approximate accurately using Monte Carlo integration. To address this issue, we propose a weak adversarial tensor neural network (WATNN) method for solving high-dimensional SFP equations, which combines a weak adversarial (WA) formulation with the tensor neural network (TNN). The WA formulation avoids the cross terms introduced by squaring the strong residual, reducing the number of integrals to be computed. When both the primal network and the adversarial test network are set as TNNs, the TNN quadrature scheme provides accurate and efficient normalization and loss computation with relatively few fixed quadrature points, avoiding the statistical sampling error that arises when using Monte Carlo integration. We demonstrate the effectiveness numerically on four high-dimensional SFP equations. For the 6D unimodal, 6D multi-modal, 10D multi-modal and 20D multi-modal problems, WATNN achieves the mean relative errors of $1.50\%$, $1.13\%$, $1.20\%$ and $7.93\%$ respectively, in only a few minutes.
\end{abstract}



\section{Introduction}
The Fokker-Planck (FP) equation describes the evolution of the probability density functions (PDFs) of dynamical systems governed by stochastic noise \cite{risken1996fokkerplanck}. It has wide applications in statistical mechanics \cite{shizgal1997quadrature}, chemistry \cite{lotstedt2006dimensional}, biology \cite{carrillo2011decision}, engineering \cite{fokou2016probabilistic} and finance \cite{gareche2013fokker}. However, the numerical methods for FP equations face three major difficulties: high dimensions, unbounded computational domains, and the normalization condition of PDFs. Grid-based methods, such as finite difference and finite element methods, suffer from the curse of dimensionality, since the number of degrees of freedom grows exponentially with the dimension. Monte Carlo methods provide an alternative approach, but reducing statistical errors can require a large number of samples in high dimensions. Moreover, reconstructing the full probability density from these samples remains challenging \cite{lei2026density}. For example, kernel density estimation also suffers from the curse of dimensionality \cite{hastie2009elements}.

Recently, deep learning methods have shown great potential for solving high-dimensional partial differential equations (PDEs) due to their universal approximation and mesh-free training. Representative approaches include physics-informed neural networks (PINNs) \cite{raissi2019physics}, the Deep Ritz method \cite{e2018deep}, weak adversarial networks (WANs) \cite{zang2020weak}, and the Deep BSDE method \cite{han2018solving}. These methods reformulate PDE problems into optimization problems and have been applied to a wide range of high-dimensional PDEs. However, applying deep learning methods to the high-dimensional stationary Fokker-Planck (SFP) equations considered in this work still involves several difficulties. For example, since the stationary equation is homogeneous, minimizing its residual alone admits the trivial zero solution, making it necessary to enforce the normalization condition. Moreover, loss computation and density normalization often involve high-dimensional integrals, which can be expensive to approximate accurately using Monte Carlo integration.

Several neural network methods have been developed for SFP equations, with different strategies for addressing the normalization difficulty. For example, the adaptive deep density approximation method in \cite{tang2022adaptive} uses KRnet to naturally satisfy the PDF constraints. The data-driven method in \cite{zhai2022deep} combines the PDE residual loss with approximate density values obtained from simulations of the associated stochastic differential equations, using the PDE constraint to reduce the reliance on accurate reference data. The physics-informed neural network method in \cite{alhussein2023physics} enforces the normalization condition by adding a weighted penalty term to the PDE residual loss. The tensor radial basis function networks in \cite{wang2025tensor} are specifically designed to satisfy the normalization condition exactly.

While the above efforts directly utilize the strong residual loss of the SFP equation, the score-based method in \cite{liu2025scorebased} instead employs a score-based PDE loss to learn an unnormalized density with free-form network architectures and enforces normalization only by post-processing. The weak generative sampler in \cite{cai2026weak} learns a transport map from a simple base distribution to the target distribution through the weak formulation of the SFP equation, thereby directly generating samples from the learned distribution.

However, most of the existing work computes the loss by Monte Carlo integration, which may require a large number of samples to accurately approximate the high-dimensional integrals involved. An alternative is to utilize the tensor neural network (TNN) \cite{wang2024tensor}, which represents a multivariate function as a finite sum of products of one-dimensional functions. This structure reduces integrals of tensor-product expressions to sums of products of one-dimensional integrals, which can be computed accurately by fixed quadrature rules. The use of this quadrature scheme for the normalization condition in \cite{liu2025scorebased,wang2025tensor} motivates its extension to loss computation. A direct approach is to combine the TNN with the strong residual loss of the SFP equation, which can be computed using the TNN quadrature scheme when the drift and diffusion coefficients admit suitable tensor-product representations. However, expanding the squared strong residual introduces cross terms that increase the number of integrals to be computed, as discussed in Remark~\ref{rem:weak_vs_strong}. The weak residual avoids this expansion and requires only first-order spatial derivatives of the networks, motivating the use of a weak adversarial (WA) formulation based on WAN \cite{zang2020weak}. Therefore, we propose the weak adversarial tensor neural network (WATNN) method for solving high-dimensional SFP equations. By setting both the primal network and the adversarial test network as TNNs, we compute the WA loss using the TNN quadrature scheme, thereby avoiding Monte Carlo integration.

The main advantages of WATNN are summarized as follows:
\begin{itemize}
    \item WATNN does not require samples from the invariant distribution for training. Unlike methods that rely on SDE simulations to generate training data, such as \cite{zhai2022deep} and the Fokker-Planck neural network (FPNN) method \cite{liu2025scorebased} adopted as the benchmark in our experiments, WATNN trains the networks directly from the weak formulation using fixed quadrature nodes. SDE simulations are used only to estimate the computational domain before training.
    \item {\emergencystretch=1em
    WATNN improves the efficiency of adversarial training by optimizing a test space rather than a single test function. We utilize the rank-one components of the test TNN to construct the test space as their linear span and determine an optimal adversarial test function within it. We further derive explicit expressions for the resulting loss and its gradient to optimize the test space.\par}
    \item WATNN utilizes the TNN quadrature scheme for accurate and efficient normalization and loss computation, avoiding Monte Carlo sampling and its statistical error. By applying high-order quadrature to the one-dimensional integrals derived from the TNN structure, the method can obtain high accuracy with fewer evaluation points than Monte Carlo integration, thereby reducing the cost of computing the high-dimensional integrals in the normalization and loss.
    \item {\emergencystretch=1em
    WATNN successfully solves the tested 6-20 dimensional SFP equations. Compared with the state-of-the-art FPNN method on the tested examples, WATNN demonstrates notable advantages in both accuracy and efficiency. For the 6D unimodal, 6D multi-modal, 10D multi-modal and 20D multi-modal problems, WATNN achieves the mean relative errors of $1.50\%$, $1.13\%$, $1.20\%$ and $7.93\%$ in only $131.28$~s, $121.50$~s, $166.94$~s and $170.87$~s, respectively.\par}
\end{itemize}

The rest of the paper is organized as follows. In Sec.~\ref{sec:pre}, we introduce the SFP equation and give a brief overview of the WAN method and TNN with its quadrature scheme, which form the basis of our proposed method. In Sec.~\ref{sec:watnn}, we develop the WATNN method, including the construction of the primal network and the test space, and the computation of the loss using the TNN quadrature scheme. In Sec.~\ref{sec:experiment}, we present four examples to demonstrate the effectiveness of WATNN. Finally, we conclude the paper in Sec.~\ref{sec:conclusion}.

\section{Preliminaries}
\label{sec:pre}
Consider the following stochastic differential equation (SDE):\begin{equation}
\label{eq:SDE}
    d\bmX_t = \bmb(\bmX_t)\,dt + \bms(\bmX_t)\,d\bmW_t,
\end{equation}
where $\bmb:\mathbb{R}^d\to\mathbb{R}^d$ is the drift coefficient, $\bms:\mathbb{R}^d\to\mathbb{R}^{d\times m}$ is the diffusion coefficient, and $\bmW_t$ is an $m$-dimensional Brownian motion. The invariant density $p$ of \eqref{eq:SDE} satisfies the stationary Fokker-Planck (SFP) equation
\begin{equation}
\label{eq:SFP}
\begin{gathered}
    \mathcal{L}p=-\sum_{i=1}^d\frac{\partial}{\partial x_i}\left(b_i p\right)(\bmx)+\frac{1}{2}\sum_{i,j=1}^d\frac{\partial^2}{\partial x_i\partial x_j}\left(D_{ij}p\right)(\bmx)=0,\\
    p(\bmx)\ge 0,\qquad \int_{\mathbb{R}^d}p(\bmx)\,d\bm{x}=1,
\end{gathered}
\end{equation}
where $\bmx = (x_1,x_2,\cdots,x_d) \in \mathbb{R}^d$, and $\bmD(\bmx)=\bms(\bmx)\bms(\bmx)^{T}$ is the diffusion matrix.

In this work, we focus on a class of SFP equations that admits a unique solution. This can be ensured under suitable structural and dissipativity conditions (see \cite{Huang2015fokkerplanck}). Since the SFP equation is posed on the entire space $\mathbb{R}^d$, our numerical method uses a bounded computational domain $\Omega\subset\mathbb{R}^d$ for approximation. We estimate $\Omega$ so that it captures almost all probability mass of $p$, namely
\begin{equation}
    \int_{\mathbb{R}^d\setminus\Omega}p(\bmx)\,d\bmx
    \le \delta,
\end{equation}
for a given small tolerance $\delta$. Consequently, we assume that the density $p(\bmx)$ decays rapidly to zero near the boundary of $\Omega$ and vanishes outside the domain.

Before introducing our WATNN method, we briefly review the WAN method and TNN with its quadrature scheme.
\subsection{Weak adversarial network}
\label{sec:wan}
A motivation for using the WAN method is to avoid directly minimizing the strong residual loss of the Fokker-Planck equation, which requires computing second-order derivatives of the network and can be computationally expensive in high dimensions, see details in Remark~\ref{rem:weak_vs_strong}. The WAN method \cite{zang2020weak} considers elliptic partial differential equations with certain boundary conditions written in the form
\begin{equation}
\label{eq:elliptic pde}
    \mathscr{A}u=f \quad \text{in } \Omega,
    \qquad
    \mathscr{B}u=g \quad \text{on } \partial\Omega,
\end{equation}
where $\mathscr{A}$ is a second-order differential operator, $\mathscr{B}$ denotes the trace operator, and $\Omega\subset\mathbb{R}^d$ is a bounded domain. We denote by $(\cdot,\cdot)_{\Omega}$ the $L^2(\Omega)$ inner product. Multiplying \eqref{eq:elliptic pde} by a test function $\varphi\in H_0^1(\Omega)$ and applying integration by parts gives the weak formulation
\begin{equation}
    a(u,\varphi) = (f,\varphi)_{\Omega},
    \qquad
    u\in H^1(\Omega),\quad \varphi\in H_0^1(\Omega).
\end{equation}
Here $a(u,\varphi)$ is obtained from the formal inner product $(\mathscr{A}u,\varphi)_{\Omega}$ after integration by parts, with the boundary terms vanishing because $\varphi=0$ on $\partial\Omega$. 

For each $u\in H^1(\Omega)$, we define the weak residual as a linear functional $\mathcal{R}[u]:H_0^1(\Omega)\to\mathbb{R}$ by
\begin{equation}
\label{eq:weak residual}
    \mathcal{R}[u](\varphi):=\langle\mathcal{R}[u],\varphi\rangle
    :=(f,\varphi)_{\Omega}-a(u,\varphi),\qquad\forall \varphi\in H_0^1(\Omega).
\end{equation}
If $u\in H^1(\Omega)$ satisfies $\mathscr{B}u=g$ and $\langle\mathcal{R}[u],\varphi\rangle=0$ for all $\varphi\in H_0^1(\Omega)$, we say that $u$ is a weak solution of \eqref{eq:elliptic pde}. The weak solution can be obtained by minimizing the following operator norm \cite{bao2020numerical}:
\begin{equation}
\label{eq:op norm}
    \|\mathcal{R}[u]\|_{\mathrm{op}}:=\sup_{\varphi\in H_0^1(\Omega),\varphi\neq 0}\frac{|\langle\mathcal{R}[u],\varphi\rangle|}{\|\varphi\|_{H^1(\Omega)}}.
\end{equation}
Thus the WAN method considers the following minimization problem:
\begin{equation}
\label{eq:classical wan}
    u^*=\mathop{\arg\min}_{u\in H^1(\Omega)}\left(\|\mathcal{R}[u]\|_{\mathrm{op}}+\beta\|\mathscr{B}u-g\|_{L^2(\partial\Omega)}\right),
\end{equation}
where $\beta>0$ is a penalty parameter. The trial solution and the test function are approximated by two neural networks,
\begin{equation}
    u(\bmx)\approx u_{\boldsymbol{\theta}}(\bmx),\qquad\varphi(\bmx)\approx \varphi_{\boldsymbol{\eta}}(\bmx),
\end{equation}
where $u_{\boldsymbol{\theta}}$ is the primal network and $\varphi_{\boldsymbol{\eta}}$ is the adversarial test network. The corresponding loss function of the WAN method is obtained by substituting the neural networks into \eqref{eq:classical wan}, and the network parameters $\boldsymbol{\theta}$ and $\boldsymbol{\eta}$ are then trained by solving the minimax problem:
\begin{equation}
\label{eq:wan_minimax}
\min_{\boldsymbol{\theta}}\max_{\boldsymbol{\eta}}\mathcal{L}_{\mathrm{WAN}}(\boldsymbol{\theta},\boldsymbol{\eta}),
\end{equation}
where the WAN loss is defined by
\begin{equation}
\label{eq:wan loss}
    \mathcal{L}_{\mathrm{WAN}}(\boldsymbol{\theta},\boldsymbol{\eta}):=\frac{|\langle\mathcal{R}[u_{\boldsymbol{\theta}}],\varphi_{\boldsymbol{\eta}}\rangle|}{\|\varphi_{\boldsymbol{\eta}}\|_{H^1(\Omega)}}+\beta\|\mathscr{B}u_{\boldsymbol{\theta}}-g\|_{L^2(\partial\Omega)}.
\end{equation}

\subsection{Tensor neural network and its quadrature scheme}
The WAN loss \eqref{eq:wan loss} introduced above still requires the computation of the weak residual $\langle\mathcal{R}[u_{\boldsymbol{\theta}}], \varphi_{\boldsymbol{\eta}}\rangle$ and the norm $\|\varphi_{\boldsymbol{\eta}}\|_{H^1(\Omega)}$. These quantities are given by integrals over the high-dimensional computational domain $\Omega$ which are usually approximated by the Monte Carlo method. However, Monte Carlo integration may require a large sample size to achieve satisfactory accuracy in high-dimensional domains. Therefore, in order to compute the weak residual and the norm efficiently in high dimensions, we consider the tensor neural network (TNN) \cite{wang2024tensor}.

TNN represents a function by a tensor-product structure. More precisely, for each dimension $i=1,\ldots,d$, we construct a one-dimensional sub-network
\begin{equation}
    \Phi_i(x_i;\boldsymbol{\theta}_i) = \left( \phi_{i,1}(x_i;\boldsymbol{\theta}_i), \phi_{i,2}(x_i;\boldsymbol{\theta}_i), \ldots, \phi_{i,R}(x_i;\boldsymbol{\theta}_i) \right)^{T},
\end{equation}
where $x_i$ is the one-dimensional input, $\boldsymbol{\theta}_i$ denotes the parameters of the $i$-th sub-network, and $R$ is the tensor rank. We simply use the multilayer perceptron (MLP) architecture \cite{hornik1989multilayer} for each subnetwork. The one-dimensional output TNN considered in this work is defined by
\begin{equation}
\label{eq:TNN}
    \Psi(\bmx;\boldsymbol{\theta}) = \sum_{r=1}^{R} \prod_{i=1}^{d} \phi_{i,r}(x_i;\boldsymbol{\theta}_i), \qquad \bmx=(x_1,\ldots,x_d)\in\Omega .
\end{equation}
Here $\boldsymbol{\theta}=\{\boldsymbol{\theta}_1,\ldots,\boldsymbol{\theta}_d\}$ denotes all parameters of the TNN. The main advantage of this structure is that the integral of a product of TNNs can be decomposed into sums of products of one-dimensional integrals. To illustrate this point, consider a general integrand written as a finite sum of products of TNNs:
\begin{equation}
\label{eq:TNN-general-product}
    F(\bmx) = \sum_{\ell=1}^{L} \prod_{m=1}^{M_{\ell}} \Psi_{\ell,m}(\bmx;\boldsymbol{\theta}^{\ell,m}),
\end{equation}
where each $\Psi_{\ell,m}$ is a TNN of the form \eqref{eq:TNN}, given by
\begin{equation}
    \Psi_{\ell,m}(\bmx;\boldsymbol{\theta}^{\ell,m}) = \sum_{r=1}^{R_{\ell,m}} \prod_{i=1}^{d} \phi^{\ell,m}_{i,r}(x_i;\boldsymbol{\theta}_i^{\ell,m}).
\end{equation}
Then each product in \eqref{eq:TNN-general-product} can be expanded as
\begin{equation}
    \prod_{m=1}^{M_{\ell}} \Psi_{\ell,m}(\bmx;\boldsymbol{\theta}^{\ell,m}) = \sum_{r_1=1}^{R_{\ell,1}} \cdots \sum_{r_{M_{\ell}}=1}^{R_{\ell,M_{\ell}}} \prod_{i=1}^{d} \left( \prod_{m=1}^{M_{\ell}} \phi^{\ell,m}_{i,r_m}(x_i;\boldsymbol{\theta}_i^{\ell,m}) \right).
\end{equation}
Therefore, the high-dimensional integral of $F$ can be written as
\begin{equation}
\label{eq:TNN-integral-splitting}
\begin{aligned}
    \int_{\Omega} F(\bmx) d\bmx &= \sum_{\ell=1}^{L} \sum_{r_1=1}^{R_{\ell,1}} \cdots \sum_{r_{M_{\ell}}=1}^{R_{\ell,M_{\ell}}} \prod_{i=1}^{d} \int_{\Omega_i} \prod_{m=1}^{M_{\ell}} \phi^{\ell,m}_{i,r_m}(x_i;\boldsymbol{\theta}_i^{\ell,m}) \,dx_i,
\end{aligned}
\end{equation}
where we assume that the domain $\Omega$ also has a tensor-product form:
\begin{equation}
    \Omega=\Omega_1\times\Omega_2\times\cdots\times\Omega_d, \qquad \Omega_i\subset\mathbb{R}.
\end{equation}
Applying one-dimensional quadrature to \eqref{eq:TNN-integral-splitting} gives the TNN quadrature scheme
\begin{equation}
\label{eq:TNN-quad}
    \int_{\Omega} F(\bmx)d\bmx \approx \sum_{\ell=1}^{L} \sum_{r_1=1}^{R_{\ell,1}} \cdots \sum_{r_{M_{\ell}}=1}^{R_{\ell,M_{\ell}}} \prod_{i=1}^{d} \left( \sum_{n=1}^{N_i} w_{i,n} \prod_{m=1}^{M_{\ell}} \phi^{\ell,m}_{i,r_m}(x_{i,n};\boldsymbol{\theta}_i^{\ell,m}) \right),
\end{equation}
where $\{x_{i,n},w_{i,n}\}_{n=1}^{N_i}$ are quadrature points and weights on $\Omega_i$. Hence, a high-dimensional integral is reduced to a finite sum of products of one-dimensional integrals. With standard one-dimensional quadrature, such as the Gauss-Legendre quadrature, \eqref{eq:TNN-quad} provides a high-accuracy approximation of the original high-dimensional integral without the need to compute the integral on the full grid.

\section{Weak adversarial tensor neural network}
\label{sec:watnn}
In this section, we present the weak adversarial tensor neural network (WATNN) method for high-dimensional SFP equations. We first estimate a bounded computational domain using simulations of the corresponding SDE and derive the WA formulation on this domain in Sec.~\ref{sec:domain_weak_formulation}. We then construct a normalized primal network and a parameterized test space using TNNs, and derive explicit expressions for the loss function and its gradient in Secs.~\ref{sec:normalized_density_network} and \ref{sec:test_space_loss}. Finally, we describe the TNN quadrature scheme for computing the loss function in Sec.~\ref{sec:tensor_quadrature} and the complete training algorithm in Sec.~\ref{sec:watnn_training}.

\subsection{Computational domain and weak adversarial formulation}
\label{sec:domain_weak_formulation}
As discussed in Sec.~\ref{sec:pre}, we approximate the invariant density $p$ on a bounded computational domain $\Omega$ containing its main probability mass. We also assume that the density $p$ decays rapidly to zero near the boundary of $\Omega$ and vanishes outside the domain. 

We estimate the computational domain by simulating the corresponding SDE dynamics of the SFP equation. More precisely, we initialize independent particles from a prescribed distribution and evolve them using the Euler-Maruyama method. With a sufficiently large terminal time and a suitably small step size, we use the particle positions at the terminal time as approximate samples from the invariant distribution \cite{mattingly2002ergodicity}. Based on the distribution of these samples, we select the product domain $\Omega=\prod_{i=1}^d[a_i,b_i]$ so that it contains a large majority of the samples. The procedure is summarized in Algorithm~\ref{alg:domain}. The domain estimation can be performed on CPUs and easily parallelized across independent SDE trajectories, requiring relatively little time under our experiment settings and is therefore excluded from the reported runtimes. Further discussion and the parameter settings are given in Sec.~\ref{sec:experiment_setup}.

\begin{algorithm}[t]
\caption{Euler-Maruyama method for estimating the computational domain}
\label{alg:domain}
\begin{algorithmic}[1]
\Require Drift $\bmb$, diffusion coefficient $\bms$, initial distribution $p_0$, particle number $Q$, terminal time $T$, step number $N$, tolerance $\delta\in(0,1)$.
\Ensure Computational domain $\Omega$.
\State Sample $\bm{x}_0^{(q)}\sim p_0$ independently for $q=1,\ldots,Q$ and set $\Delta t=T/N$.
\For{$n=0,\ldots,N-1$}
    \For{$q=1,\ldots,Q$}
        \State Sample $\Delta\bm{W}_n^{(q)}\sim\mathcal{N}(\bm{0},\Delta t\boldsymbol{I}_{m\times m})$ independently.
        \State $\bm{x}_{n+1}^{(q)}\gets\bm{x}_n^{(q)}+\bmb(\bm{x}_n^{(q)})\Delta t+\bms(\bm{x}_n^{(q)})\Delta\bm{W}_n^{(q)}$.
    \EndFor
\EndFor
\State Select intervals $[a_i,b_i]$, $i=1,\ldots,d$, such that $\prod_{i=1}^d[a_i,b_i]$ contains at least $Q(1-\delta)$ terminal samples.
\State \Return $\Omega=\prod_{i=1}^d[a_i,b_i]$.
\end{algorithmic}
\end{algorithm}

Next we derive the weak formulation of the SFP equation on the computational domain $\Omega$. Multiplying the SFP equation \eqref{eq:SFP} by a test function $\varphi\in H^1(\Omega)$ and applying integration by parts gives
\begin{align}
\label{eq:int by parts}
    (\mathcal{L}p,\varphi)_\Omega
    &=-\frac{1}{2}\sum_{i,j=1}^d\int_\Omega\frac{\partial(D_{ij}p)}{\partial x_j}\frac{\partial\varphi}{\partial x_i}\,d\bm{x}+\sum_{i=1}^d\int_\Omega b_i p\frac{\partial\varphi}{\partial x_i}\,d\bm{x} \notag\\
    &\quad-\sum_{i=1}^d\int_{\partial\Omega}b_i p\varphi n_i\,dS+\frac{1}{2}\sum_{i,j=1}^d\int_{\partial\Omega}\frac{\partial(D_{ij}p)}{\partial x_j}\varphi n_i\,dS,
\end{align}
where $\bm{n}=(n_1,\ldots,n_d)$ is the outward unit normal vector on $\partial\Omega$. The boundary integrals involving $b_i p$ and $\partial(D_{ij}p)/\partial x_j$ are negligible due to the decay assumption of the density $p$. Thus the test functions need not vanish on $\partial\Omega$ and the test space is $H^1(\Omega)$ rather than $H_0^1(\Omega)$ as in Sec.~\ref{sec:wan}. Omitting these boundary integrals yields the weak residual functional $\mathcal{R}[p]:H^1(\Omega)\to\mathbb{R}$,
\begin{equation}
\label{eq:sfp_weak_residual}
    \langle\mathcal{R}[p],\varphi\rangle:=-\frac{1}{2}\sum_{i,j=1}^d\int_\Omega\frac{\partial(D_{ij}p)}{\partial x_j}\frac{\partial\varphi}{\partial x_i}\,d\bm{x}+\sum_{i=1}^d\int_\Omega b_i p\frac{\partial\varphi}{\partial x_i}\,d\bm{x}.
\end{equation}
We therefore seek the density $p$ satisfying
\begin{equation}
\label{eq:sfp_truncated_weak}
    \langle\mathcal{R}[p],\varphi\rangle=0 \quad\forall\varphi\in H^1(\Omega).
\end{equation}
 
We approximate the density $p$ by a nonnegative, normalized primal network $p_{\boldsymbol{\theta}}$ and consider a finite-dimensional test space $V_{\boldsymbol{\eta}}\subset H^1(\Omega)$ spanned by basis functions parameterized by $\boldsymbol{\eta}$. The specific constructions of the primal network and the test space are described below. Following the WAN formulation in Sec.~\ref{sec:wan}, we define the loss function
\begin{equation}
\label{eq:watnn_sfp_loss}
    \mathcal{L}_{\mathrm{WAN}}(\boldsymbol{\theta},\boldsymbol{\eta}):=\sup_{\substack{\varphi\in V_{\boldsymbol{\eta}}\\\varphi\neq0}}\frac{|\langle\mathcal{R}[p_{\boldsymbol{\theta}}],\varphi\rangle|}{\|\varphi\|_{H^1(\Omega)}}.
\end{equation}
For fixed $\boldsymbol{\theta}$ and $\boldsymbol{\eta}$, the supremum corresponds to selecting an optimal test function from $V_{\boldsymbol{\eta}}$. We solve the minimax problem
\begin{equation}
\label{eq:watnn_sfp_minimax}
    \min_{\boldsymbol{\theta}}\max_{\boldsymbol{\eta}}\mathcal{L}_{\mathrm{WAN}}(\boldsymbol{\theta},\boldsymbol{\eta}).
\end{equation}
The primal network is trained to minimize this loss, while the test space is optimized to maximize it. Computing the loss requires high-dimensional integrals for the weak residual and the $H^1$ norm, which motivates the use of TNN.

\subsection{Normalized primal TNN}
\label{sec:normalized_density_network}
To construct the normalized primal network $p_{\boldsymbol{\theta}}$, we normalize the one-dimensional factors of a TNN and combine the resulting tensor-product densities with trainable weights. For the domain $\Omega=\prod_{i=1}^d[a_i,b_i]$ and each coordinate $x_i\in[a_i,b_i]$, we use an MLP sub-network $\widetilde{\bm{p}}_i(x_i;\boldsymbol{\theta}_i)$ with $r_p$ outputs,
\begin{equation}
\label{eq:density_factor_network}
    \widetilde{\bm{p}}_i(x_i;\boldsymbol{\theta}_i)=\left(\widetilde p_{i,1}(x_i;\boldsymbol{\theta}_i),\ldots,\widetilde p_{i,r_p}(x_i;\boldsymbol{\theta}_i)\right)^T, \qquad i=1,\ldots,d.
\end{equation}
We apply the softplus activation $\sigma(z)=\ln(1+e^z)$ to the output layer to ensure positivity. To eliminate the boundary terms in \eqref{eq:int by parts}, we require the primal network and its first spatial derivatives to vanish on $\partial\Omega$. For the domain $\Omega=\prod_{i=1}^d[a_i,b_i]$, we multiply each output $\widetilde p_{i,r}$ by 
\begin{equation}
    \max\left(\left(1-\frac{(2x_i-a_i-b_i)^2}{(b_i-a_i)^2}\right)^3,0\right),
    \qquad i=1,\ldots,d,
\end{equation}
before normalization to ensure the boundary condition. We retain the notation $\widetilde p_{i,r}$ for the resulting outputs.

We normalize each output over its corresponding interval by setting
\begin{equation}
\label{eq:normalized_density_factors}
    Z_{i,r}(\boldsymbol{\theta}_i)=\int_{a_i}^{b_i}\widetilde p_{i,r}(x_i;\boldsymbol{\theta}_i)\,dx_i, \qquad p_{i,r}(x_i;\boldsymbol{\theta}_i)=\frac{\widetilde p_{i,r}(x_i;\boldsymbol{\theta}_i)}{Z_{i,r}(\boldsymbol{\theta}_i)}, \qquad r=1,\ldots,r_p.
\end{equation}
Each $p_{i,r}$ is thus a probability density on $[a_i,b_i]$, and its tensor product over the coordinates is a probability density on $\Omega$. We define the normalized primal network by
\begin{equation}
\label{eq:normalized_density_tnn}
    p_{\boldsymbol{\theta}}(\bm{x})=\sum_{r=1}^{r_p}\alpha_r(\boldsymbol{\beta})\prod_{i=1}^d p_{i,r}(x_i;\boldsymbol{\theta}_i), \qquad \alpha_r(\boldsymbol{\beta})=\frac{e^{\beta_r}}{\sum_{s=1}^{r_p}e^{\beta_s}},
\end{equation}
where $\boldsymbol{\beta}\in\mathbb{R}^{r_p}$ and $\boldsymbol{\theta}=(\boldsymbol{\theta}_1,\ldots,\boldsymbol{\theta}_d,\boldsymbol{\beta})$ denotes the trainable parameters of the primal network $p_{\boldsymbol{\theta}}$. It is easy to verify that $p_{\boldsymbol{\theta}}$ is nonnegative and satisfies the normalization condition
\begin{equation}
\label{eq:density_tnn_unit_mass}
    \int_\Omega p_{\boldsymbol{\theta}}(\bm{x})\,d\bm{x}=\sum_{r=1}^{r_p}\alpha_r(\boldsymbol{\beta})\prod_{i=1}^d\int_{a_i}^{b_i}p_{i,r}(x_i;\boldsymbol{\theta}_i)\,dx_i=\sum_{r=1}^{r_p}\alpha_r(\boldsymbol{\beta})=1.
\end{equation}
In implementation, each normalizing constant $Z_{i,r}$ in \eqref{eq:normalized_density_factors} is approximated by the composite Gauss-Legendre quadrature on $[a_i,b_i]$, as described in Sec.~\ref{sec:tensor_quadrature}. 

\subsection{Test space and the loss function}
\label{sec:test_space_loss}
We construct the test space $V_{\boldsymbol{\eta}}$ from the rank-one components of a test TNN $\varphi_{\boldsymbol{\eta}}$. For the domain $\Omega=\prod_{i=1}^d[a_i,b_i]$ and each coordinate $x_i\in[a_i,b_i]$, the test sub-network has $r_\varphi$ outputs,
\begin{equation}
\label{eq:test_factor_network}
    \widetilde{\bm{\varphi}}_i(x_i;\boldsymbol{\eta}_i)=\left(\widetilde\varphi_{i,1}(x_i;\boldsymbol{\eta}_i),\ldots,\widetilde\varphi_{i,r_\varphi}(x_i;\boldsymbol{\eta}_i)\right)^T, \qquad i=1,\ldots,d,
\end{equation}
where $r_\varphi$ is also the tensor rank. To control the scales of the rank-one components, we normalize each output in $L^2([a_i,b_i])$,
\begin{equation}
\label{eq:normalized_test_factors}
\begin{aligned}
    s_{i,r}(\boldsymbol{\eta}_i)&=\left(\int_{a_i}^{b_i}\widetilde\varphi_{i,r}(x_i;\boldsymbol{\eta}_i)^2\,dx_i\right)^{1/2},\\
    \varphi_{i,r}(x_i;\boldsymbol{\eta}_i)&=\frac{\widetilde\varphi_{i,r}(x_i;\boldsymbol{\eta}_i)}{s_{i,r}(\boldsymbol{\eta}_i)}, \qquad r=1,\ldots,r_\varphi.
\end{aligned}
\end{equation}
For fixed $\boldsymbol{\eta}$, this normalization only rescales each rank-one component and therefore preserves its linear span.

With $\boldsymbol{\eta}=(\boldsymbol{\eta}_1,\ldots,\boldsymbol{\eta}_d)$, we define the rank-one functions and their linear span by
\begin{equation}
\label{eq:tnn_test_space}
\begin{aligned}
    \psi_r(\bm{x};\boldsymbol{\eta})&=\prod_{i=1}^d\varphi_{i,r}(x_i;\boldsymbol{\eta}_i), \qquad r=1,\ldots,r_\varphi,\\
    V_{\boldsymbol{\eta}}&=\operatorname{span}\{\psi_1,\ldots,\psi_{r_\varphi}\}\subset H^1(\Omega).
\end{aligned}
\end{equation}
Thus, the supremum in \eqref{eq:watnn_sfp_loss} ranges over linear combinations
\begin{equation}
\label{eq:test_linear_combination}
    \varphi_{\boldsymbol{c},\boldsymbol{\eta}}(\bm{x})=\sum_{r=1}^{r_\varphi}c_r\psi_r(\bm{x};\boldsymbol{\eta}), \qquad \boldsymbol{c}\in\mathbb{R}^{r_\varphi}.
\end{equation}

Define the residual vector $\boldsymbol{e}=\boldsymbol{e}(\boldsymbol{\theta},\boldsymbol{\eta})$ and the Gram matrix $G=G(\boldsymbol{\eta})$ by
\begin{equation}
\label{eq:test_residual_gram}
\begin{aligned}
    e_s&=\langle\mathcal{R}[p_{\boldsymbol{\theta}}],\psi_s\rangle,\\
    G_{rs}&=(\psi_r,\psi_s)_{H^1(\Omega)}
    =\int_\Omega\left(\psi_r\psi_s+\nabla\psi_r\cdot\nabla\psi_s\right)\,d\bm{x}, \qquad r,s=1,\ldots,r_\varphi.
\end{aligned}
\end{equation}
The linearity of the weak residual in the test function gives
\begin{equation}
\label{eq:test_coefficient_residual_norm}
    \langle\mathcal{R}[p_{\boldsymbol{\theta}}],\varphi_{\boldsymbol{c},\boldsymbol{\eta}}\rangle=\boldsymbol{c}^T\boldsymbol{e}, \qquad \|\varphi_{\boldsymbol{c},\boldsymbol{\eta}}\|_{H^1(\Omega)}^2=\boldsymbol{c}^TG\boldsymbol{c}.
\end{equation}
If the functions $\{\psi_r\}_{r=1}^{r_\varphi}$ are linearly independent, $G$ is symmetric positive definite. Applying the Cauchy-Schwarz inequality in the inner product induced by $G$ then yields
\begin{equation}
\label{eq:wan_gram_loss}
    \mathcal{L}_{\mathrm{WAN}}(\boldsymbol{\theta},\boldsymbol{\eta})=\sup_{\boldsymbol{c}\neq\boldsymbol{0}}\frac{|\boldsymbol{c}^T\boldsymbol{e}|}{\sqrt{\boldsymbol{c}^TG\boldsymbol{c}}}=\sqrt{\boldsymbol{e}^TG^{-1}\boldsymbol{e}}.
\end{equation}
To compute the loss and its gradient, we solve
\begin{equation}
\label{eq:wan_energy_solve}
    G\boldsymbol c=\boldsymbol e, \qquad \mathcal L_{\mathrm{WAN}}=\sqrt{\boldsymbol e^T\boldsymbol c}.
\end{equation}
For any scalar parameter $\xi$ in $\boldsymbol\theta$ or $\boldsymbol\eta$, differentiating $G\boldsymbol c=\boldsymbol e$ gives $G\,\partial\boldsymbol c/\partial\xi=\partial\boldsymbol e/\partial\xi-(\partial G/\partial\xi)\boldsymbol c$. Substitution into \eqref{eq:wan_energy_solve} yields
\begin{equation}
\label{eq:wan_loss_gradient}
    \frac{\partial\mathcal L_{\mathrm{WAN}}}{\partial\xi}=\frac{1}{2\mathcal L_{\mathrm{WAN}}}\left(2\boldsymbol c^T\frac{\partial\boldsymbol e}{\partial\xi}-\boldsymbol c^T\frac{\partial G}{\partial\xi}\boldsymbol c\right), \qquad \mathcal L_{\mathrm{WAN}}>0.
\end{equation}
The derivatives of $\boldsymbol e$ and $G$ in \eqref{eq:wan_loss_gradient} with respect to the network parameters are computed by automatic differentiation. Thus \eqref{eq:wan_loss_gradient} gives the loss gradient without differentiating through the linear solver, so the solver need not support automatic differentiation.

\subsection{Computation of the loss function using the TNN quadrature scheme}
\label{sec:tensor_quadrature}
{\emergencystretch=1em
To compute the loss function \eqref{eq:wan_gram_loss}, we apply the TNN quadrature scheme \eqref{eq:TNN-quad} to the residual vector and the Gram matrix. For each interval $[a_i,b_i]$, let $\{x_{i,n},w_{i,n}\}_{n=1}^{N_i}$ be the nodes and weights of the composite Gauss-Legendre quadrature, and define
\begin{equation}
\label{eq:watnn_one_dim_quadrature}
    \mathcal I_i[g]:=\sum_{n=1}^{N_i}w_{i,n}g(x_{i,n})\approx\int_{a_i}^{b_i}g(x_i)\,dx_i.
\end{equation}
}
The normalizing integrals in \eqref{eq:normalized_density_factors} and \eqref{eq:normalized_test_factors} are approximated by
\begin{equation}
\label{eq:factor_normalization_quadrature}
    Z_{i,r}\approx\mathcal I_i[\widetilde p_{i,r}], \qquad s_{i,r}\approx\left(\mathcal I_i[\widetilde\varphi_{i,r}^{\,2}]\right)^{1/2}.
\end{equation}

The tensor-product structure of $\psi_r$ gives the following quadrature approximation of the $H^1$ Gram matrix entries $G_{rs}$ for $r,s=1,\ldots,r_\varphi$:
\begin{equation}
\label{eq:test_gram_quadrature}
    G_{rs}\approx\prod_{i=1}^d\mathcal I_i[\varphi_{i,r}\varphi_{i,s}]
    +\sum_{j=1}^d\mathcal I_j\left[\frac{\partial\varphi_{j,r}}{\partial x_j}\frac{\partial\varphi_{j,s}}{\partial x_j}\right]\prod_{i\ne j}\mathcal I_i[\varphi_{i,r}\varphi_{i,s}].
\end{equation}
Thus, all entries of the Gram matrix $G$ are approximated by one-dimensional inner products of the test factors and their first derivatives.

Finally, from \eqref{eq:sfp_weak_residual}, computing each entry $e_s=\langle\mathcal{R}[p_{\boldsymbol{\theta}}],\psi_s\rangle$ of the residual vector $\boldsymbol{e}$ requires approximating the integrals
\begin{equation}
\label{eq:residual_integral_forms}
    \int_\Omega b_i(\bm x)p_{\boldsymbol\theta}(\bm x)\frac{\partial\psi_s}{\partial x_i}(\bm x)\,d\bm x
    \quad\text{and}\quad
    \int_\Omega\frac{\partial}{\partial x_j}\left(D_{ij}(\bm x)p_{\boldsymbol\theta}(\bm x)\right)\frac{\partial\psi_s}{\partial x_i}(\bm x)\,d\bm x.
\end{equation}
Suppose the coefficients admit TNN approximations
\begin{equation}
\label{eq:coefficient_tensor_representations}
    b_i(\bm x)\approx\sum_{r=1}^{R_i^b}\prod_{k=1}^d b_{i,r,k}(x_k), \qquad D_{ij}(\bm x)\approx\sum_{r=1}^{R_{ij}^D}\prod_{k=1}^d d_{ij,r,k}(x_k).
\end{equation}
Using \eqref{eq:normalized_density_tnn} and \eqref{eq:tnn_test_space}, each drift term is approximated by
\begin{equation}
\label{eq:residual_tensor_quadrature}
\begin{aligned}
    \int_\Omega b_i(\bm x)p_{\boldsymbol\theta}(\bm x)\frac{\partial\psi_s}{\partial x_i}(\bm x)\,d\bm x\approx\sum_{r=1}^{r_p}\alpha_r\sum_{t=1}^{R_i^b}\mathcal I_i\left[b_{i,t,i}p_{i,r}\frac{\partial\varphi_{i,s}}{\partial x_i}\right]\prod_{k\ne i}\mathcal I_k[b_{i,t,k}p_{k,r}\varphi_{k,s}].
\end{aligned}
\end{equation}
And the diffusion term is approximated by
\begin{equation}
\label{eq:diffusion_tensor_quadrature}
\begin{aligned}
    &\int_\Omega\frac{\partial}{\partial x_j}\left(D_{ij}(\bm x)p_{\boldsymbol\theta}(\bm x)\right)\frac{\partial\psi_s}{\partial x_i}(\bm x)\,d\bm x\\
    \approx &
    \begin{cases}
    \displaystyle\sum_{r=1}^{r_p}\alpha_r\sum_{t=1}^{R_{ij}^D}\mathcal I_i\left[\frac{\partial(d_{ii,t,i}p_{i,r})}{\partial x_i}\frac{\partial\varphi_{i,s}}{\partial x_i}\right]\prod_{k\ne i}\mathcal I_k[d_{ii,t,k}p_{k,r}\varphi_{k,s}], & i=j,\\
    \begin{aligned}
        &\sum_{r=1}^{r_p}\alpha_r\sum_{t=1}^{R_{ij}^D}\mathcal I_i\left[d_{ij,t,i}p_{i,r}\frac{\partial\varphi_{i,s}}{\partial x_i}\right]\mathcal I_j\left[\frac{\partial(d_{ij,t,j}p_{j,r})}{\partial x_j}\varphi_{j,s}\right]\\
        &\qquad{}\times\prod_{k\notin\{i,j\}}\mathcal I_k[d_{ij,t,k}p_{k,r}\varphi_{k,s}],
    \end{aligned}
    & i\ne j.
    \end{cases}
\end{aligned}
\end{equation}
All terms are thus computed through one-dimensional quadrature, with only first-order spatial derivatives computed by automatic differentiation.
\begin{remark}
\label{rem:weak_vs_strong}
For the SFP equation, consider the strong residual loss 
\begin{equation}
\label{eq:strong_residual_loss}
\begin{aligned}
    \int_\Omega|\mathcal Lp_{\boldsymbol\theta}|^2\,d\bm x
    =&\sum_{i,k=1}^d\int_\Omega
    \frac{\partial(b_i p_{\boldsymbol\theta})}{\partial x_i}
    \frac{\partial(b_k p_{\boldsymbol\theta})}{\partial x_k}\,d\bm x
    -\sum_{i,j,k=1}^d\int_\Omega
    \frac{\partial(b_i p_{\boldsymbol\theta})}{\partial x_i}
    \frac{\partial^2(D_{jk}p_{\boldsymbol\theta})}{\partial x_j\partial x_k}\,d\bm x\\
    &+\frac14\sum_{i,j,k,l=1}^d\int_\Omega
    \frac{\partial^2(D_{ij}p_{\boldsymbol\theta})}{\partial x_i\partial x_j}
    \frac{\partial^2(D_{kl}p_{\boldsymbol\theta})}{\partial x_k\partial x_l}\,d\bm x.
\end{aligned}
\end{equation}
When the coefficients admit TNN approximations, each integral can be computed using the TNN quadrature scheme \eqref{eq:TNN-quad}. However, the strong residual loss \eqref{eq:strong_residual_loss} contains $d^2+d^3+d^4$ integrals, compared with $d+d^2$ integrals in the weak residual \eqref{eq:sfp_weak_residual}. Moreover, the weak residual has only first-order spatial derivatives of the networks. We therefore construct the loss using the weak residual rather than the strong residual.
\end{remark}

\subsection{Training algorithm}
\label{sec:watnn_training}
The complete WATNN procedure is summarized in Algorithm~\ref{alg:watnn1}. The computational domain $\Omega=\prod_{i=1}^d[a_i,b_i]$ is estimated from SDE trajectories by Algorithm~\ref{alg:domain}, and composite Gauss-Legendre quadrature nodes and weights are generated on each interval $[a_i,b_i]$, $i=1,\ldots,d$. At each training iteration, the test space is first updated for $K_\eta$ optimizer steps to maximize the loss with the primal network fixed, after which the primal network is updated for $K_\theta$ optimizer steps to minimize the loss with the test space fixed. After training is finished, $p_{\boldsymbol\theta}$ is used as the numerical approximation to the invariant density.

\begin{algorithm}[t]
\caption{WATNN for solving the SFP equation}
\label{alg:watnn1}
\begin{algorithmic}[1]
\item[\textbf{Input:}] Computational domain $\Omega=\prod_{i=1}^d[a_i,b_i]$ estimated by Algorithm~\ref{alg:domain}, quadrature sizes $N_i$, training iterations $K$, update numbers $K_\theta$ and $K_\eta$, and learning rates.
\item[\textbf{Output:}] Approximation $p_{\boldsymbol\theta}$ to the invariant density of the SFP equation.
\State Initialize the primal network $p_{\boldsymbol\theta}$ and test space $V_{\boldsymbol{\eta}}$.
\State For each interval $[a_i,b_i]$, generate composite Gauss-Legendre quadrature nodes and weights $\{x_{i,n},w_{i,n}\}_{n=1}^{N_i}$.
\For{$k=1,\ldots,K$}
    \For{$j=1,\ldots,K_\eta$}
        \State Compute residual vector $\boldsymbol e$ and Gram matrix $G$ using the TNN quadrature scheme as described in Sec.~\ref{sec:tensor_quadrature}.
        \State Compute the loss by \eqref{eq:wan_energy_solve} and its gradient with respect to $\boldsymbol\eta$ by \eqref{eq:wan_loss_gradient}.
        \State Update $\boldsymbol\eta$ by one optimizer step to maximize the loss, keeping $\boldsymbol\theta$ fixed.
    \EndFor
    \State Compute $G$ at the updated $\boldsymbol{\eta}$ using the TNN quadrature scheme.
    \For{$j=1,\ldots,K_\theta$}
        \State Compute $\boldsymbol e$ using the TNN quadrature scheme.
        \State Compute the loss by \eqref{eq:wan_energy_solve} and its gradient with respect to $\boldsymbol\theta$ by \eqref{eq:wan_loss_gradient}.
        \State Update $\boldsymbol\theta$ by one optimizer step to minimize the loss, keeping $\boldsymbol\eta$ fixed.
    \EndFor
\EndFor
\State \Return $p_{\boldsymbol\theta}$.
\end{algorithmic}
\end{algorithm}

\section{Numerical Experiments}
\label{sec:experiment}
In this section, we evaluate the performance of WATNN on four benchmark problems for high-dimensional SFP equations. Specifically, we consider the 6D unimodal, 6D multi-modal, and 10D multi-modal problems from \cite{wang2025tensor}, together with the 20D multi-modal problem constructed in this work. These examples span dimensions from 6 to 20 and combine different invariant density structures with nonlinear interactions among coordinates, providing challenging tests for high-dimensional SFP solvers.

\subsection{Experiment setup}
\label{sec:experiment_setup}
Following a construction similar to that in \cite{wang2025tensor}, we consider benchmark SFP equations with known invariant densities. Specifically, we take \eqref{eq:SFP} with the drift
\begin{equation}
   \label{eq:drift}
    \bmb(\bmx) = -\frac{1}{2}\bmD(\bmx)\nabla H(\bmx)
    +\frac{1}{2}\nabla\cdot\bmD(\bmx),
\end{equation}
where $H:\mathbb{R}^d\to\mathbb{R}$ is the potential function and the matrix divergence $(\nabla\cdot\bmD)_i=\sum_{j=1}^d\partial D_{ij}/\partial x_j$. The corresponding invariant density $p^*(\bmx)$ is then given by
\begin{equation}
    p^*(\bmx) = \frac{1}{Z}e^{-H(\bmx)}, \qquad Z = \int_{\mathbb{R}^d} e^{-H(\bmx)} d\bm{x}.
\end{equation}
Although WATNN is applicable to more general drift fields, all problems considered here use the form in \eqref{eq:drift} to evaluate the performance of WATNN using known invariant densities.

We compare WATNN with the Fokker-Planck neural network (FPNN) method \cite{liu2025scorebased}, considering both TNN-based and MLP-based implementations of FPNN. FPNN learns an unnormalized density by minimizing a score-based PDE loss, whose expectation with respect to the invariant distribution is approximated using training samples. After training is finished, the normalization constant is computed using the TNN quadrature scheme for TNN-based FPNN and Monte Carlo integration for MLP-based FPNN. We train FPNN with Adam \cite{kingma2014adam} at a learning rate of $0.01$. For the 6D and 10D problems, the training dataset contains $20k$ points, with a batch size of $2k$. For the more challenging 20D problem, we increase the training dataset size to $200k$ and the batch size to $20k$. For MLP-based FPNN, we also increase the number of Monte Carlo samples used to estimate the normalization constant from $20k$ to $1000k$ for the 20D problem. These increases are made in an effort to reduce sampling errors and obtain more accurate FPNN solutions.

The prediction accuracy is measured by the mean relative error (MRE),
\begin{equation}
    \label{eq:mre}
    \frac{1}{n} \sum_{i=1}^n\left|\frac{p^*\left(\boldsymbol{x}_i\right)-p_{\boldsymbol{\theta}}\left(\boldsymbol{x}_i\right)}{p^*\left(\boldsymbol{x}_i\right)}\right|, \quad \boldsymbol{x}_i \in \mathcal{D}_{\text{test}}=\left\{\boldsymbol{x} \in \mathbb{R}^d \mid p^*(\boldsymbol{x})>\varepsilon\right\} .
\end{equation}
Here $p_{\boldsymbol{\theta}}$ denotes the neural network approximation and $p^*$ is the reference density. The test dataset $\mathcal{D}_{\text{test}}$ is generated by using the gradient-ascent method described in \cite{liu2025scorebased}, with a learning rate of $1\times10^{-3}$, $10^4$ initial points and a threshold $\varepsilon$ given in Table~\ref{tab:experiment_parameters}.

In all experiments, the sub-networks of the primal TNN $p_{\boldsymbol{\theta}}$ and the test TNN $\varphi_{\boldsymbol{\eta}}$ use the MLP architecture. And we construct the test space from the rank-one components of the test TNN $\varphi_{\boldsymbol{\eta}}$ as described in Sec.~\ref{sec:test_space_loss}. The hidden layers use tanh activation in the primal sub-networks and sinc activation in the test sub-networks. A softplus activation is applied to the output layer of each primal sub-network to enforce positivity. We also adopt the power embedding in \cite{liu2025scorebased}. For each dimension $i$, the scalar input $x_i$ is transformed into $[x_i,x_i^2,\ldots,x_i^m]$ before being fed into the corresponding sub-network to enhance the representation of nonlinear dependence.

We use the notation $[m,a_1,\ldots,a_\ell,r]$ to describe the architectures of both the primal TNN and the test TNN used to construct the test space. Here, $m$ is the power embedding order, $a_1,\ldots,a_\ell$ are the widths of the $\ell$ hidden layers in each sub-network, and $r$ is the tensor rank, equal to the number of outputs of each sub-network. For example, $[5,64,64,64]$ in Table~\ref{tab:watnn_6d_unimodal} denotes a rank-$64$ TNN with power embedding order $5$, whose sub-networks have two hidden layers both of width $64$ and $64$ outputs.

Unless otherwise stated, both the primal TNN and the test space are trained with Adam. Across all four problems, the learning rates follow polynomial decay schedules, decreasing from $1\times10^{-3}$ to $1\times10^{-5}$ for the primal TNN and from $1\times10^{-2}$ to $1\times10^{-4}$ for the test TNN.

We select the computational domains based on Algorithm~\ref{alg:domain}. In this algorithm, we initialize $2\times 10^4$ particles uniformly on $[-0.5,0.5]^d$, simulate the particles up to time $T=20$ with step size $\Delta t=5\times 10^{-3}$, and set the tolerance to $\delta=0.001$. For the TNN quadrature scheme, each coordinate interval is divided into $M$ equal subintervals, with 16 Gauss-Legendre points on each subinterval, giving $N=16M$ points per coordinate. The different settings for each problem are shown in Table~\ref{tab:experiment_parameters}, where the small thresholds are chosen to evaluate the performance of WATNN in relatively low density regions.

\begin{table}[htbp]
\centering
\caption{Settings for the four benchmark problems. Here $\varepsilon$ denotes the density threshold for selecting test points with $p^*(\bm{x})>\varepsilon$. Each coordinate interval of domain $\Omega$ is divided into $M$ equal subintervals, giving $N=16M$ quadrature points per coordinate.}
\label{tab:experiment_parameters}
\begin{tabular}{lllll}
\toprule
Problem & Domain $\Omega$ & Threshold $\varepsilon$ & $M$ & $N$ \\
\midrule
6D unimodal & $[-1.2,1.2]^6$ & $10^{-5}$ & 6 & 96 \\
6D multi-modal & $[-4,4]^6$ & $10^{-5}$ & 10 & 160 \\
10D multi-modal & $[-2,2]^{10}$ & $10^{-8}$ & 6 & 96 \\
20D multi-modal & $[-1.2,1.2]^{20}$ & $10^{-8}$ & 6 & 96 \\
\bottomrule
\end{tabular}
\end{table}

Finally, WATNN computations are performed in double precision to maintain the numerical accuracy of the TNN quadrature scheme. All WATNN and FPNN training runs are performed on a single NVIDIA A100 GPU with 80~GB of memory. Computational domain estimation for WATNN is excluded from the reported runtimes (see the discussion in Sec.~\ref{sec:domain_weak_formulation}), as is training data generation for FPNN for the same reason. The runtime for WATNN is measured from the start of program execution, whereas that for FPNN includes only the training time.

\subsection{6D Unimodal problem}
We first consider the drift \eqref{eq:drift} with the potential function
\begin{equation}
\label{eq:6d_unimodal_potential}
    H(\bm{x})
    =
    3\left[
    (x_1^4-x_2)^2+2x_2^2
    +(x_3^4-x_4)^2+2x_4^2
    +(x_5^4-x_6)^2+2x_6^2
    \right],
\end{equation}
and the diffusion $\bmD(\bmx) = 2\boldsymbol{I}_{6\times6}$ on the selected domain $\Omega=[-1.2,1.2]^6$. The couplings between pairs of variables in \eqref{eq:6d_unimodal_potential} lead to a high-order nonlinear drift, making this example challenging for neural network training. 

Tables~\ref{tab:watnn_6d_unimodal} and \ref{tab:fpnn_6d_unimodal} demonstrate the advantage of WATNN in both prediction accuracy and computational efficiency. In particular, with the $[8,64,128]$ architecture, WATNN achieves its best MRE of $1.50\%$ in $131.28$ s. The TNN-based FPNN yields substantially larger errors at higher computational cost, while the MLP-based FPNN also remains less accurate at comparable runtimes.

\begin{figure}[htbp]
    \centering
    \begin{subfigure}[b]{0.32\textwidth}
        \centering
        \includegraphics[width=\linewidth]{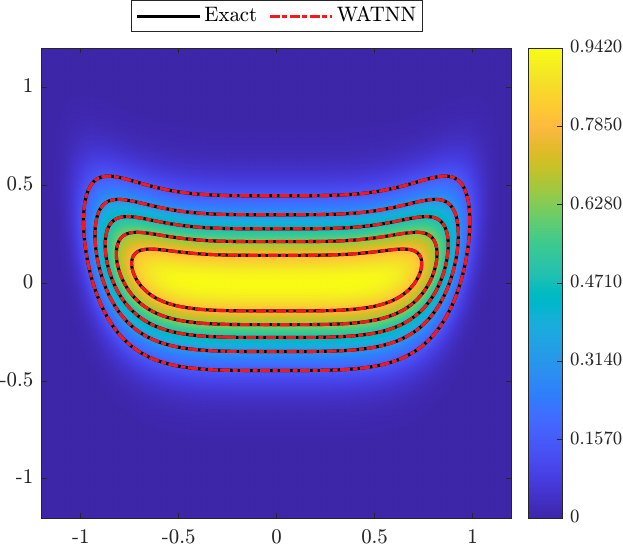}
        \caption{$(x_1,x_2,0,0,0,0)$}
    \end{subfigure}
   % \qquad
    \begin{subfigure}[b]{0.32\textwidth}
        \centering
        \includegraphics[width=\linewidth]{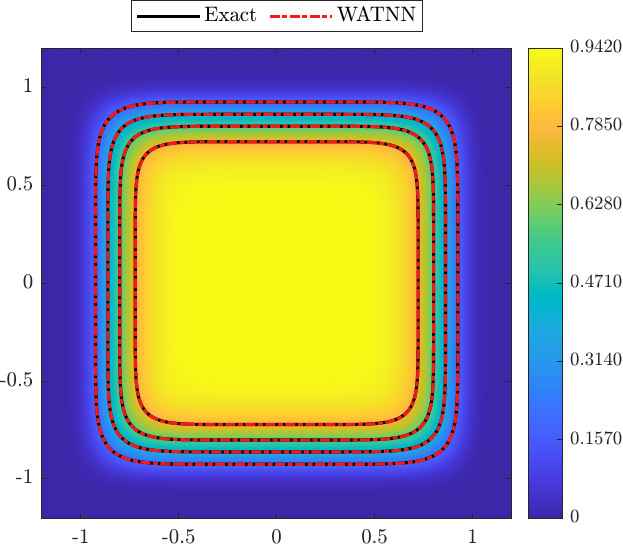}
        \caption{$(x_1,0,x_3,0,0,0)$}
    \end{subfigure} 
    \begin{subfigure}[b]{0.32\textwidth}
        \centering
        \includegraphics[width=\linewidth]{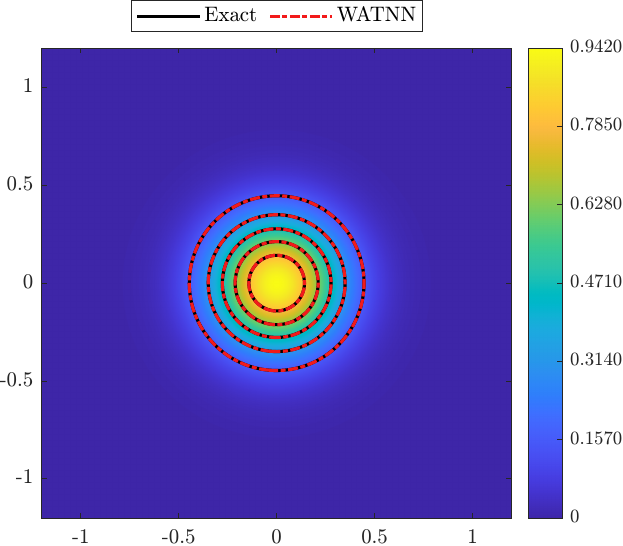}
        \caption{ $(0,x_2,0,x_4,0,0)$}
    \end{subfigure}
    %\qquad
    \caption{Contour plots of projections of the exact and WATNN solutions for the 6D unimodal problem.}
    \label{fig:6d_unimodal}
\end{figure}

We present contour plots of different projections of the exact solution and the best WATNN result obtained with the $[8,64,128]$ architecture in Fig.~\ref{fig:6d_unimodal}, where two coordinates are varied while the remaining coordinates are fixed at zero. The close agreement of the contours indicates that WATNN accurately captures the shape of the unimodal density. 

\begin{table}[htbp]
\centering
\caption{Results of WATNN on the 6D unimodal problem. $K_\theta$ and $K_\eta$ denote the numbers of primal network and test space updates per step, respectively.}
\label{tab:watnn_6d_unimodal}
\begin{tabular}{lllllll}
\toprule
 Primal TNN & Test TNN & Steps & $K_\theta$ & $K_\eta$ & MRE & Time (s)\\ \midrule
 $[m, \text{hidden layers}, r]$ &  &  &  &  &  &\\
 $[1, 64, 64, 64]$ & $[1, 64, 64, 64]$ & $20k$ & 1& 1  & 7.67\% & 34.22 \\
 $[3, 64, 64, 64]$ & $[3, 64, 64, 64]$ & $20k$ & 1& 1 & 4.96\% & 34.66 \\
 $[5, 64, 64, 64]$ & $[5, 64, 64, 64]$ & $20k$ & 1& 1 & 5.59\% & 34.73 \\
 $[7, 64, 64, 64]$ & $[7, 64, 64, 64]$ & $20k$ & 1& 1 & 5.71\% & 34.31 \\
 $[8, 64, 128]$ & $[8, 64, 128]$ & $20k$ & 1 & 1  & \textbf{2.29\%} & 40.00 \\
 \midrule
 $[1, 64, 64, 64]$ & $[1, 64, 64, 64]$ & $80k$ & 1& 1  & 4.76\% & 108.22 \\
 $[3, 64, 64, 64]$ & $[3, 64, 64, 64]$ & $80k$ & 1& 1 & 3.66\% & 110.02 \\
 $[5, 64, 64, 64]$ &  $[5, 64, 64, 64]$ & $80k$ & 1& 1 & 3.79\% & 110.12 \\
 $[7, 64, 64, 64]$ &  $[7, 64, 64, 64]$ & $80k$ & 1& 1 & 3.40\% & 108.81 \\
 $[8, 64, 128]$ & $[8, 64, 128]$ & $80k$ & 1& 1  & \textbf{1.50\%} & 131.28 \\
\bottomrule
\end{tabular}
\end{table} 

\begin{table}[htbp]
\centering
\caption{Results of FPNN with TNN/MLP architectures on the 6D unimodal problem.}
\label{tab:fpnn_6d_unimodal}
\begin{tabular}{lllll}
\toprule
Model & Layers & Steps & MRE & Time (s) \\
\midrule
 & $[m,\text{hidden layers},r]$ & & & \\
 & $[1,64,64,64]$ & $2k$ & 30.61\% & 149.38 \\
 & $[3,64,64,64]$ & $2k$ & 10.45\% & 160.08 \\
TNN & $[5,64,64,64]$ & $2k$ & 9.34\% & 191.84 \\
 & $[7,64,64,64]$ & $2k$ & 9.95\% & 198.99 \\
\cmidrule(l){2-5}
 & $[8,64,128]$ & $2k$ & 7.69\% & 160.64 \\
 & $[8,64,128]$ & $8.5k$ & 4.69\% & 665.98 \\
\midrule
 & $[d,\text{hidden layers}]$ & & & \\
 & $[6,20,20,20]$ & $5k$ & 17.51\% & 83.95 \\
MLP & $[6,20,20,20]$ & $10k$ & 14.36\% & 167.24 \\
\cmidrule(l){2-5}
 & $[6,64,64,64]$ & $5k$ & 13.40\% & 92.85 \\
 & $[6,64,64,64]$ & $10k$ & 8.63\% & 179.30 \\
\bottomrule
\end{tabular}
\end{table}

\subsection{6D Multi-modal problem}
We next consider a 6D multi-modal problem with drift \eqref{eq:drift} and diffusion matrix $\bmD(\bmx)=2\boldsymbol{I}_{6\times6}$ on the selected domain $\Omega=[-4,4]^6$. The potential function is given by
\begin{equation}
\label{eq:6d_multimodal_potential}
\begin{aligned}
    H(\bm{x})
    ={}&2\left[
    x_1^2+x_2^2+x_3^2
    +0.5\left(x_1x_2+x_1x_3+x_2x_3\right)
    \right] \\
    &-\ln\left(x_1^2+0.02\right)
    -\ln\left(x_2^2+0.02\right) \\
    &+0.5\left[
    x_4^2+x_5^2+x_6^2
    +0.2\left(x_4x_5+x_4x_6+x_5x_6\right)
    \right].
\end{aligned}
\end{equation}
The coordinates interact within the groups $(x_1,x_2,x_3)$ and $(x_4,x_5,x_6)$. The four isolated local minima of the potential give rise to four isolated modes of the invariant density, requiring the approximation to resolve both the separation of the high-density regions and the coordinate interactions. 

The results in Tables~\ref{tab:watnn_6d_multimodal} and \ref{tab:fpnn_6d_multimodal} show that WATNN retains its advantage when the density has multiple modes. The best result is obtained with the $[5,64,64,64]$ architecture, which achieves an MRE of $1.13\%$ in only $121.50$ s. For comparison, the MLP-based FPNN obtains its lowest MRE of $13.88\%$ in $124.11$ s, which remains substantially higher than that of WATNN at a comparable runtime.

\begin{figure}[htbp]
    \centering
    \begin{subfigure}[b]{0.32\textwidth}
        \centering
        \includegraphics[width=\linewidth]{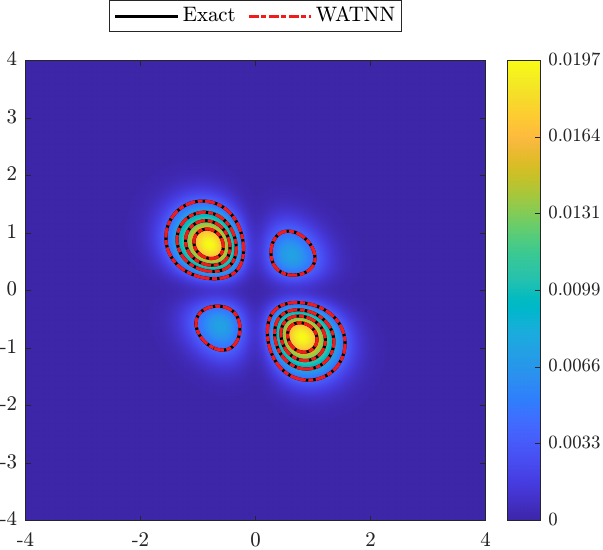}
        \caption{$(x_1,x_2,0,0,0,0)$}
    \end{subfigure}
   % \qquad
    \begin{subfigure}[b]{0.32\textwidth}
        \centering
        \includegraphics[width=\linewidth]{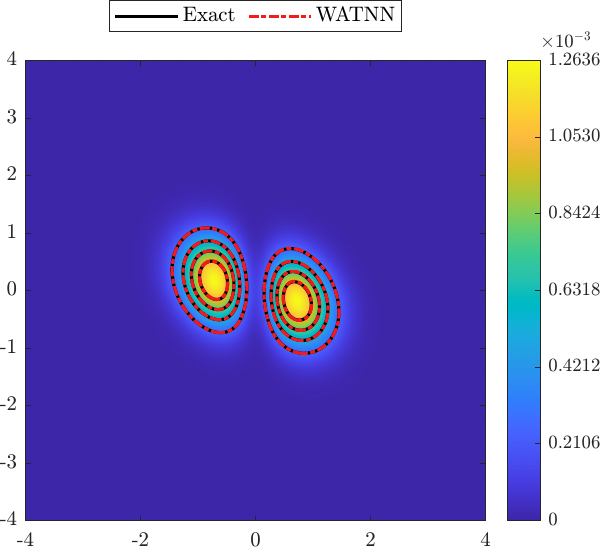}
        \caption{$(x_1,0,x_3,0,0,0)$}
    \end{subfigure} 
    \begin{subfigure}[b]{0.32\textwidth}
        \centering
        \includegraphics[width=\linewidth]{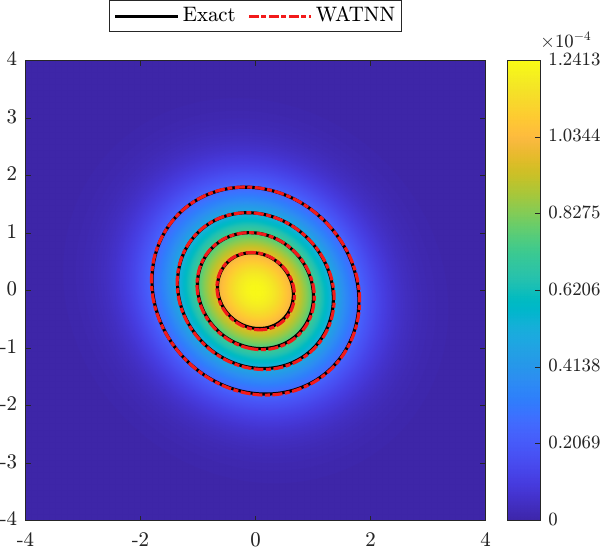}
        \caption{ $(0,0,0,x_4,x_5,0)$}
    \end{subfigure}
    %\qquad
    \caption{Contour plots of projections of the exact and WATNN solutions for the 6D multi-modal problem.}
    \label{fig:6d_multimodal}
\end{figure}

We further present the contour plots of different projections of the exact solution and the best WATNN result in Fig.~\ref{fig:6d_multimodal}, where the WATNN result closely matches the exact solution, indicating that WATNN successfully reproduces the complex multi-modal structure.

\begin{table}[htbp]
\centering
\caption{Results of WATNN on the 6D multi-modal problem. $K_\theta$ and $K_\eta$ denote the numbers of primal network and test space updates per step, respectively.}
\label{tab:watnn_6d_multimodal}
\begin{tabular}{lllllll}
\toprule
Primal TNN & Test TNN & Steps & $K_\theta$ & $K_\eta$ & MRE & Time (s) \\
\midrule
$[m, \text{hidden layers}, r]$ &  &  &  &  &  &\\
$[1,64,64,64]$ & $[1,64,64,64]$ & $20k$ & 1 & 1 & 17.31\% & 36.64 \\
$[3,64,64,64]$ & $[3,64,64,64]$ & $20k$ & 1 & 1 & \textbf{4.07\%} & 37.93 \\
$[5,64,64,64]$ & $[5,64,64,64]$ & $20k$ & 1 & 1 & 4.52\% & 37.71 \\
$[7,64,64,64]$ & $[7,64,64,64]$ & $20k$ & 1 & 1 & 7.43\% & 37.94 \\
$[8,64,128]$ & $[8,64,128]$ & $20k$ & 1 & 1 & 9.93\% & 46.72 \\
\midrule
$[1,64,64,64]$ & $[1,64,64,64]$ & $80k$ & 1 & 1 & 15.52\% & 118.93 \\
$[3,64,64,64]$ & $[3,64,64,64]$ & $80k$ & 1 & 1 & 1.64\% & 122.53 \\
$[5,64,64,64]$ & $[5,64,64,64]$ & $80k$ & 1 & 1 & \textbf{1.13\%} & 121.50 \\
$[7,64,64,64]$ & $[7,64,64,64]$ & $80k$ & 1 & 1 & 1.51\% & 122.65 \\
$[8,64,128]$ & $[8,64,128]$ & $80k$ & 1 & 1 & 4.08\% & 156.61 \\
\bottomrule
\end{tabular}
\end{table}

\begin{table}[htbp]
\centering
\caption{Results of FPNN with TNN/MLP architectures on the 6D multi-modal problem.}
\label{tab:fpnn_6d_multimodal}
\begin{tabular}{lllll}
\toprule
Model & Layers & Steps & MRE & Time (s)\\ \midrule
& $[m, \text{hidden layers}, r]$ &  &  &  \\
 & $[1, 64, 64, 64]$  & $5k$ & 34.55\% & 381.09\\
TNN & $[3, 64, 64, 64]$  & $5k$ & 96.48\%  & 401.19\\
 & $[5, 64, 64, 64]$  & $5k$ & 23.01\%& 429.22\\
 & $[8, 64, 128]$  & $5k$ & 26.53\% & 388.92\\ \midrule
& $[d, \text{hidden layers}]$ &  &  &   \\
 & $[6, 20, 20, 20]$ & $5k$ & 15.61\%& 76.82\\
MLP & $[6, 32, 32, 32]$ & $5k$ & 16.02\%& 77.41\\
 & $[6, 64, 64]$  & $5k$ & 16.29\% & 64.11\\ 
  & $[6, 64, 64]$  & $10k$ & 13.88\% & 124.11\\
\bottomrule
\end{tabular}
\end{table}

\subsection{10D Multi-modal problem}
We now turn to a higher-dimensional example with a multi-modal invariant density in $\mathbb{R}^{10}$. Consider the drift \eqref{eq:drift} with the potential function
\begin{equation}
\label{eq:10d_multimodal_potential}
\begin{aligned}
    H(\bm{x})
    ={}&2.5\left[
    x_1^2+x_2^2+x_3^2
    +0.1\left(x_1x_2+x_1x_3+x_2x_3\right)
    \right] \\
    &+2\left[
    x_4^2+x_5^2+x_6^2
    +0.2\left(x_4x_5+x_4x_6+x_5x_6\right)
    \right] \\
    &+3\left(x_7^2+x_8^2-0.01x_7x_8\right)
    +3\left(x_9^2+x_{10}^2-0.01x_9x_{10}\right) \\
    &-\ln\left(2x_9^2+0.02\right),
\end{aligned}
\end{equation}
and the diffusion $\bmD(\bmx)=2\boldsymbol{I}_{10\times10}$ on the selected domain $\Omega=[-2,2]^{10}$. Here, the coordinate interactions occur within four groups, namely $(x_1,x_2,x_3)$, $(x_4,x_5,x_6)$, $(x_7,x_8)$, and $(x_9,x_{10})$. The potential has two isolated local minima, corresponding to two modes of the invariant density. The higher dimensionality poses an additional challenge for approximating the multi-modal density.

Tables~\ref{tab:watnn_10d_multimodal} and \ref{tab:fpnn_10d_multimodal} show that WATNN maintains satisfactory accuracy and computational efficiency on this more challenging 10D problem, while FPNN performs much worse than on the 6D problems. FPNN yields much larger MREs than WATNN, and the TNN-based FPNN with the $[3,64,64,64]$ architecture even becomes numerically unstable, returning NaN after only $100$ steps. This instability may be related to the fact that FPNN does not explicitly enforce density normalization during training and can produce excessively large or small values with the multiplicative structure of TNN.

\begin{figure}[htbp]
    \centering
    \begin{subfigure}[b]{0.32\textwidth}
        \centering
        \includegraphics[width=\linewidth]{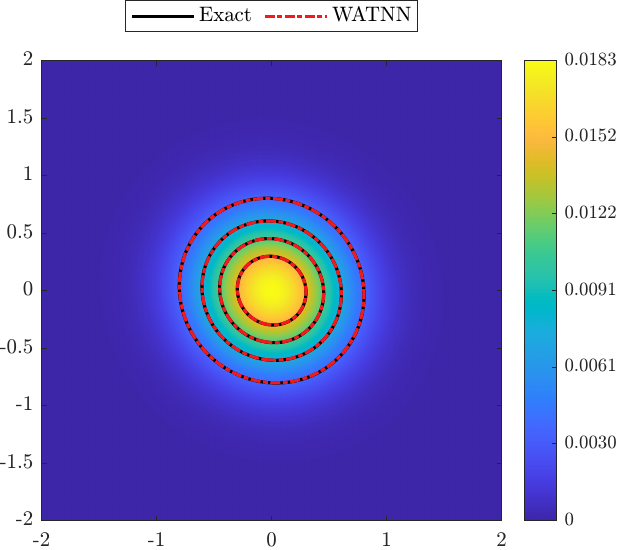}
        \caption{$(x_1,x_2,0,0,0,0,0,0,0,0)$}
    \end{subfigure}
   % \qquad
    \begin{subfigure}[b]{0.32\textwidth}
        \centering
        \includegraphics[width=\linewidth]{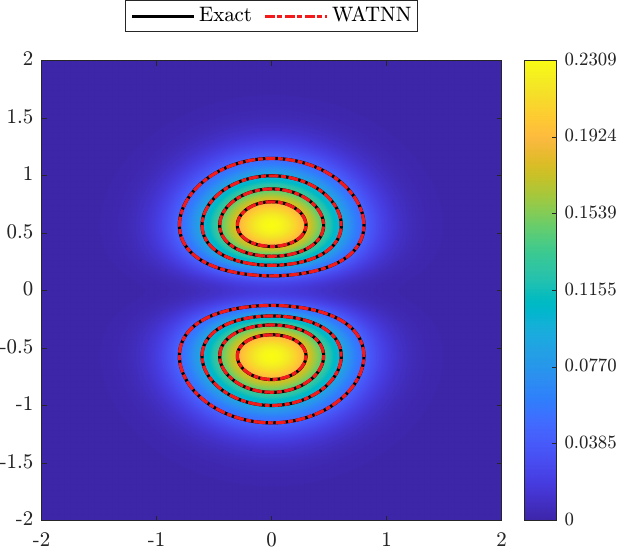}
        \caption{$(x_1,0,0,0,0,0,0,0,x_9,0)$}
    \end{subfigure} 
    \begin{subfigure}[b]{0.32\textwidth}
        \centering
        \includegraphics[width=\linewidth]{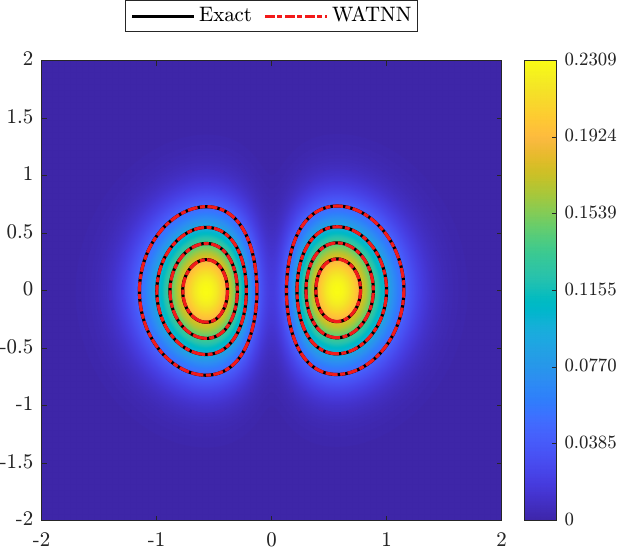}
        \caption{ $(0,0,0,0,0,0,0,0,x_9,x_{10})$}
    \end{subfigure}
    %\qquad
    \caption{Contour plots of projections of the exact and WATNN solutions for the 10D multi-modal problem.}
    \label{fig:10d_multimodal}
\end{figure}

In contrast, with the $[7,64,64,64]$ architecture for both the primal and test networks, WATNN achieves an MRE of $2.40\%$ in only $64.50$ s, substantially outperforming FPNN in terms of accuracy while requiring less computational time. After further training, the MRE is reduced to $1.20\%$ in $166.94$ s, demonstrating its potential for higher-accuracy approximation. Moreover, the contour plots of different projections of the exact solution and the best WATNN result in Fig.~\ref{fig:10d_multimodal} indicate that WATNN successfully reproduces the multi-modal structure of the exact solution.

\begin{table}[htbp]
\centering
\caption{Results of WATNN on the 10D multi-modal problem. $K_\theta$ and $K_\eta$ denote the numbers of primal network and test space updates per step, respectively.}
\label{tab:watnn_10d_multimodal}
\begin{tabular}{lllllll}
\toprule
Primal TNN & Test TNN & Steps & $K_\theta$ & $K_\eta$ & MRE & Time (s) \\
\midrule
$[m, \text{hidden layers}, r]$ &  &  &  &  &  &\\
$[1,64,64,64]$ & $[7,64,64,64]$ & $20k$ & 1 & 1 & 5.17\% & 47.17 \\
$[3,64,64,64]$ & $[7,64,64,64]$ & $20k$ & 1 & 1 & 3.28\% & 48.16 \\
$[5,64,64,64]$ & $[7,64,64,64]$ & $20k$ & 1 & 1 & 4.21\% & 48.02 \\
$[7,64,64,64]$ & $[7,64,64,64]$ & $20k$ & 1 & 1 & 6.04\% & 46.93 \\
$[7,64,64,64]$ & $[7,64,64,64]$ & $20k$ & 1 & 2 & \textbf{2.40\%} & 64.50 \\
$[8,64,128]$ & $[7,64,64,64]$ & $20k$ & 1 & 1 & 3.78\% & 54.90 \\
\midrule
$[1,64,64,64]$ & $[7,64,64,64]$ & $80k$ & 1 & 1 & 2.71\% & 153.37 \\
$[3,64,64,64]$ & $[7,64,64,64]$ & $80k$ & 1 & 1 & 1.60\% & 156.58 \\
$[5,64,64,64]$ & $[7,64,64,64]$ & $80k$ & 1 & 1 & 1.96\% & 156.24 \\
$[7,64,64,64]$ & $[7,64,64,64]$ & $80k$ & 1 & 1 & 1.97\% & 154.58 \\
$[7,64,64,64]$ & $[7,64,64,64]$ & $60k$ & 1 & 2 & \textbf{1.20\%} & 166.94 \\
$[8,64,128]$ & $[7,64,64,64]$ & $80k$ & 1 & 1 & 1.76\% & 180.67 \\
\bottomrule
\end{tabular}
\end{table}

\begin{table}[htbp]
\centering
\caption{Results of FPNN with TNN/MLP architectures on the 10D multi-modal problem.}
\label{tab:fpnn_10d_multimodal}
\begin{tabular}{lllll}
\toprule
Model & Layers & Steps & MRE & Time (s)\\ \midrule
& $[m, \text{hidden layers}, r]$ &  &  &  \\
 & $[1, 64, 64, 64]$  & $0.5k$ & 32.00\% & 103.96\\
TNN & $[1, 64, 64, 64]$  & $1k$ & 34.07\%  & 203.07\\ 
& $[3, 64, 64, 64]$  & $0.1k$ & NaN & 22.38\\  \midrule
& $[d, \text{hidden layers}]$ &  &  &   \\
 & $[10,64,64]$ & $5k$ & 32.40\%& 101.92\\
MLP & $[10,20,20,20,20]$  & $5k$ & 26.62\% & 135.72\\ 
 & $[10,64,64,64,64]$ & $5k$ & 20.74\%& 142.22\\
\bottomrule
\end{tabular}
\end{table}

\subsection{20D Multi-modal problem}
{\emergencystretch=1em
Finally, we assess the performance of WATNN in much higher dimensions using a 20D problem with drift \eqref{eq:drift} and diffusion $\bmD(\bmx)=2\boldsymbol{I}_{20\times20}$ on the selected domain $\Omega=[-1.2,1.2]^{20}$. Let $\bm{x}^{(k)}=\allowbreak(x_{4k-3},\allowbreak x_{4k-2},\allowbreak x_{4k-1},\allowbreak x_{4k})$, $k=1,\ldots,5$, and define
\begin{equation}
  Q_c(\bm{z})=\sum_{i=1}^{4}z_i^4+0.1\cdot c(z_1^2z_2^2+z_2^2z_3^2+z_3^2z_4^2+z_4^2z_1^2).
\end{equation}
}
The potential takes the form
\begin{equation}
\label{eq:20d_multimodal_potential}
H(\bm{x})
=
5\sum_{k=1}^{5}Q_{k}\!\left(\bm{x}^{(k)}\right)
-\ln\left(2x_{13}^2+0.02\right)
-\ln\left(2x_{17}^2+0.02\right).
\end{equation}
Here, the coordinates interact within five groups, introducing nontrivial couplings among different variables. The potential has four isolated local minima, corresponding to four distinct modes of the invariant density. The combination of high dimensionality, nonlinear interactions, and multiple modes poses significant difficulties for the numerical approximation of the invariant density. We follow the training setup in Sec.~\ref{sec:experiment_setup}, except that the test space is trained with Lion~\cite{chen2023symbolic} instead of Adam.

FPNN performs much worse on this problem, even with the larger training dataset and batch size specified in Sec.~\ref{sec:experiment_setup}. As reported in Table~\ref{tab:fpnn_20d_multimodal}, the TNN-based FPNN returns NaN after $50$ steps, while the MLP-based FPNN yields MREs of $95.88\%$ and $73.02\%$. These results further illustrate the difficulty of this problem and expose the limitations of FPNN in solving it under the tested settings.

\begin{figure}[htbp]
    \centering
    \begin{subfigure}[b]{0.32\textwidth}
        \centering
        \includegraphics[width=\linewidth]{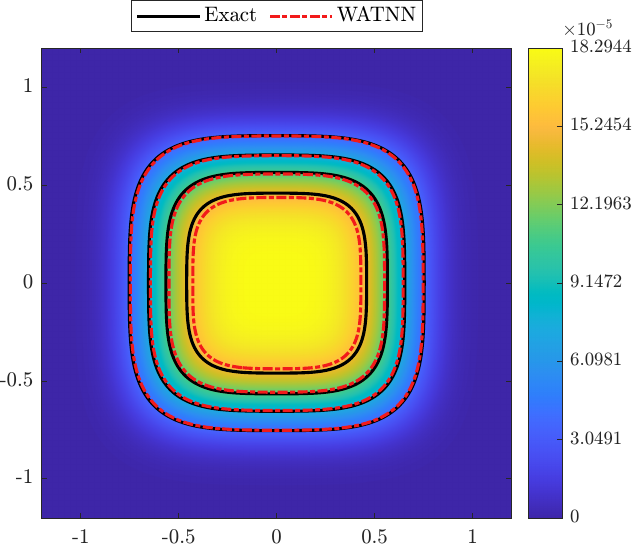}
        \caption{$(x_1,x_2;\,\mathbf{0})$}
    \end{subfigure}
   % \qquad
    \begin{subfigure}[b]{0.32\textwidth}
        \centering
        \includegraphics[width=\linewidth]{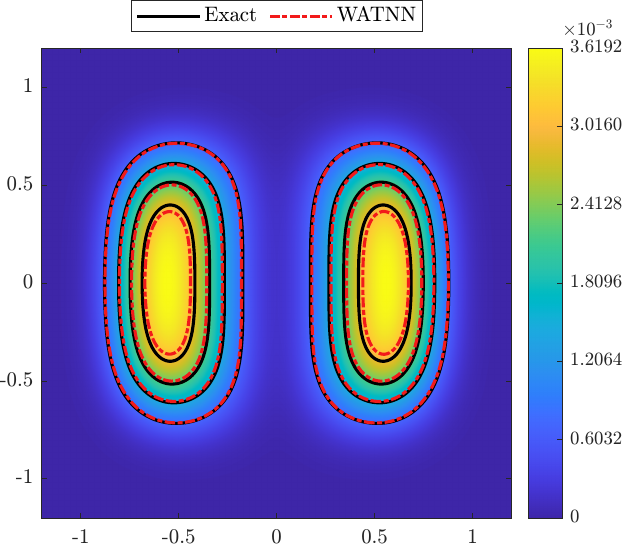}
        \caption{$(x_{13},x_{14};\,\mathbf{0})$}
    \end{subfigure} 
    \begin{subfigure}[b]{0.32\textwidth}
        \centering
        \includegraphics[width=\linewidth]{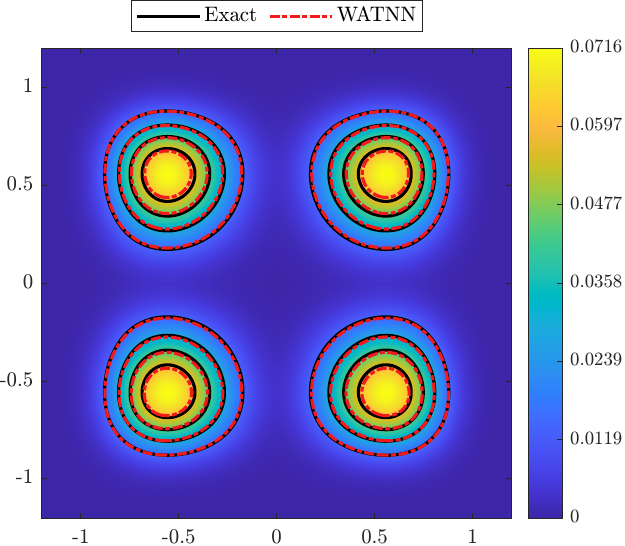}
        \caption{$(x_{13},x_{17};\,\mathbf{0})$}
    \end{subfigure}
    %\qquad
    \caption{Contour plots of projections of the exact and WATNN solutions for the 20D multi-modal problem, where $\mathbf{0}$ indicates that the remaining coordinates are fixed at zero.}
    \label{fig:20d_multimodal}
\end{figure}

WATNN still achieves satisfactory performance on this problem. As shown in Table~\ref{tab:watnn_20d_multimodal}, it achieves an MRE of $7.93\%$ in only $170.87$~s after $20k$ steps using the $[7,64,64,64]$ architecture. Further training to $80k$ steps reduces the MRE to $4.31\%$, with a total runtime of $610.88$~s. The projections in Fig.~\ref{fig:20d_multimodal} further support the ability of WATNN to reproduce the stationary density in this challenging high-dimensional setting, including its multi-modal structure in the $(x_{13},x_{17})$ plane. Moreover, WATNN also accurately approximates regions of low density, as shown in the $(x_1,x_2)$ plane, where the density is on the order of $10^{-4}$.

\begin{table}[htbp]
\centering
\caption{Results of WATNN on the 20D multi-modal problem. $K_\theta$ and $K_\eta$ denote the numbers of primal network and test space updates per step, respectively.}
\label{tab:watnn_20d_multimodal}
\begin{tabular}{lllllll}
\toprule
Primal TNN & Test TNN & Steps & $K_\theta$ & $K_\eta$ & MRE & Time (s) \\
\midrule
$[1,64,64,64]$ & $[1,64,64,64]$ & $20k$ & 3 & 1 & 42.70\% & 149.20 \\
$[3,64,64,64]$ & $[3,64,64,64]$ & $20k$ & 1 & 3 & 8.62\% & 173.64 \\
$[5,64,64,64]$ & $[5,64,64,64]$ & $20k$ & 1 & 3 & 11.79\% & 171.48 \\
$[7,64,64,64]$ & $[7,64,64,64]$ & $20k$ & 1 & 3 & \textbf{7.93\%} & 170.87 \\
$[8,64,128]$ & $[7,64,64,64]$ & $20k$ & 1 & 3 & 23.12\% & 237.01 \\
\midrule
$[1,64,64,64]$ & $[1,64,64,64]$ & $80k$ & 3 & 1 & 23.11\% & 514.19 \\
$[3,64,64,64]$ & $[3,64,64,64]$ & $80k$ & 1 & 3 & 4.90\% & 622.44 \\
$[5,64,64,64]$ & $[5,64,64,64]$ & $80k$ & 1 & 3 & 7.34\% & 617.36 \\
$[7,64,64,64]$ & $[7,64,64,64]$ & $80k$ & 1 & 3 & \textbf{4.31\%} & 610.88 \\
$[8,64,128]$ & $[7,64,64,64]$ & $80k$ & 1 & 3 & 9.70\% & 877.88 \\
\bottomrule
\end{tabular}
\end{table}

\begin{table}[htbp]
\centering
\caption{Results of FPNN with TNN/MLP architectures on the 20D multi-modal problem.}
\label{tab:fpnn_20d_multimodal}
\begin{tabular}{lllll}
\toprule
Model & Layers & Steps & MRE & Time (s)\\ \midrule
& $[m, \text{hidden layers}, r]$ &  &  &  \\
TNN & $[1, 64, 64, 64]$  & $0.05k$ & NaN & 42.37\\  \midrule
& $[d, \text{hidden layers}]$ &  &  &   \\
MLP & $[20,20,20,20,20]$  & $5k$ & 95.88\% & 644.05\\ 
 & $[20,64,64,64,64]$ & $5k$ & 73.02\%& 656.63\\
\bottomrule
\end{tabular}
\end{table}

\section{Conclusion}
\label{sec:conclusion}
In this work, we have developed the weak adversarial tensor neural network (WATNN) method for solving high-dimensional SFP equations. The method applies the WAN framework to the weak formulation of the SFP equation, with both the primal and test networks represented by TNNs. To improve the efficiency of adversarial training, we construct and optimize a test space spanned by the rank-one components of the test TNN. We also utilize the TNN quadrature scheme to enforce the normalization condition and compute the loss function with relatively few fixed quadrature points, avoiding the statistical sampling error that arises when using Monte Carlo integration. Numerical experiments on challenging 6-20 dimensional SFP problems demonstrate the strong potential of WATNN and its advantages over the existing methods in both accuracy and efficiency. 

However, when the coefficients or nonlinear terms of the equation do not admit suitable tensor-product representations, the TNN quadrature scheme cannot be directly applied, limiting the applicability of WATNN to more general equations. Although WATNN training does not rely on the probabilistic interpretation of the SFP equation, the current domain estimation still relies on such an interpretation through SDE simulations. This dependence may also restrict its application to more general equations on unbounded domains. Combining WATNN with adaptive Hermite approximations on unbounded domains using the techniques in \cite{xia2021efficient} may eliminate this dependence and allow extensions to equations without a probabilistic interpretation. Exploring these directions will be left for future work.

\section*{Acknowledgments}
This work is supported by the National Natural Science Foundation of China under the Grant No. 12325112, and the High-performance Computing Platform of Peking University.
\FloatBarrier 
\bibliographystyle{siamplain}
\bibliography{references.bib}
\end{document}